\documentclass[review]{elsarticle}

\usepackage{hyperref}
\usepackage{amsmath}
\usepackage{amsfonts}
\usepackage{amssymb}
\usepackage{makeidx}
\usepackage{graphicx}
\usepackage{color}
\usepackage[ruled]{algorithm2e}
\usepackage[margin=1.5in]{geometry}

\newtheorem{theorem}{Theorem}[section]
\newtheorem{lemma}[theorem]{Lemma}

\newtheorem{definition}[theorem]{Definition}
\newtheorem{example}[theorem]{Example}
\numberwithin{equation}{section}

\newcommand{\w}{\omega}
\newcommand{\s}{\mathbf{s}}
\newcommand{\x}{\mathbf{x}}
\newcommand{\uin}{\tiny \mbox{in}}
\newcommand{\bu}{\mathbf{u}}

\newcommand{\A}{\mathcal{A}}
\newcommand{\F}{\mathcal{F}}
\newcommand{\Sc}{\mathcal{S}}
\newcommand{\bv}{\mathbf{v}}
\newcommand{\vertiii}[1]{{\left\vert\kern-0.25ex\left\vert\kern-0.25ex\left\vert #1 
		\right\vert\kern-0.25ex\right\vert\kern-0.25ex\right\vert}}

\def\eop{{\ \vrule height 6pt width 6pt depth 0pt}}

\begin{document}

\begin{frontmatter}

\title{A  discrete-ordinate weak Galerkin  method  for  radiative transfer equation}

\author{Maneesh Kumar Singh  }
\ead{maneesh-kumar.singh@imperial.ac.uk}
\address[focal]{Department of Mathematics, Imperial College London,  SW7 2AZ, London, UK}

\begin{abstract}
This research article discusses a  numerical solution of the radiative transfer equation based on the weak Galerkin finite element method. We discretize the angular variable by means of the discrete-ordinate method. Then the resulting semi-discrete hyperbolic system is approximated using the weak Galerkin method.  The stability result for the proposed numerical method is devised.  A \emph{priori} error analysis is established under the suitable norm. In order to examine the theoretical results, numerical experiments are carried out. 
\end{abstract}

\begin{keyword}
 Radiative transfer  equation \sep Weak Galerkin method  \sep Discrete-ordinate method \sep Stability results \sep Error analysis
\MSC[2020]  65N12 \sep 65N15 \sep 65N30 \sep 65R20
\end{keyword}

\end{frontmatter}


\section{Introduction}

The propagation, absorption, and scattering of particles within the participating media are simulated using the radiative transfer equation (RTE).  The application of RTE takes place in numerous scientific fields, including atmospheric physics, biomedical optics, heat transfer, neutron transport, optical imaging, radiative transport, etc. 
The RTE is an integro-differential equation that depends on a set of optical parameters (index of refraction, absorption, scattering, and scattering function) that describe the medium, further information in \cite{ishimaru1978wave}. Schuster first devised the radiative transfer theory to analyse light propagation in foggy atmospheres \cite{schuster1905radiation}.  Chandrasekhar \cite{chandrasekhar}  contributed greatly to the theory of radiative transfer
and one of his main contributions was the creation of the vector equation of transfer for polarised light, focusing on Rayleigh scattering. In recent years,  this theory has received significant attention for its applications in combustion system \cite{modest2016radiative}, remote sensing \cite{tsang1984radiative, karam1988electromagnetic}, optical tomography \cite{klose2002optical_1,klose2002optical_2},  photo-bioreactors \cite{pilon2011radiation}  and solar energy harvesting \cite{mahian2013review, benoit2016review}. 
A comprehensive overview of the analytical and numerical solution of RTE in different participating mediums can be found in \cite{modest2013radiative}  and references therein.

\subsection{The model problem}
This research explores a stationary radiative transfer equation that is modulated by 
\begin{equation}\label{SRte1}
\left\{
\begin{array}{ll}
\displaystyle \s \cdot \nabla u(\x,\s)+\sigma_{t}(\x)u(\x,\s)-\sigma_{s}(\x) \mathcal{K}u(\x,\s) =f(\x,\s),
\quad (\x,\s) \in \Omega \times \Sc,\\[4pt]
u(\x,\s)=u_{\uin}(\x,\s),\,\,\, (\x,\s) \in \Gamma^{-},
\end{array}
\right.
\end{equation}

where $u$ is the radiative intensity  for all $\x$ in the  domain $\Omega$,   an open bounded polytopal domain in $\mathbb{R}^{d}(d=1,2,3)$  with a smooth boundary $\partial \Omega$ and  for all directions $\s$ of the unit sphere
$\Sc$ in $\mathbb{R}^{d}$.  
The incoming and outgoing boundaries are given by
\[
\Gamma^{\pm}=\{(\x,\s):\x \in \partial \Omega_{\s}^{\pm},\,\s \in \Sc\},
\]
where the spatial boundaries $\Omega_{\s}^{\pm}$ are
\[
\partial \Omega_{\s}^{-}=\{\x \in \partial \Omega: \s \cdot \mathbf{n}(\x)<0\},\quad \partial \Omega_{\s}^{+}=\{\x \in \partial \Omega: \s \cdot \mathbf{n}(\x) \geq 0\},
\]
where $\mathbf{n}(\x)$
the unit outward normal for $\x \in \partial \Omega$.  

The physical meaning of the medium parameters that determine the RTE is explained in the following few lines. 
\begin{itemize}
    \item \textbf{Absorption:} The absolute amount of absorption is directly proportional to the magnitude of the incident energy and the distance the beam travels through the medium.
    \item \textbf{Scattering:}  or  ``out-scattering'' (away from the direction $\s$ under consideration), is very similar to absorption, \emph{i.e.}, a part of the incoming intensity is removed from the direction of propagation, $\s'$.
    \item \textbf{Total attenuation:} The total attenuation of the intensity in a pencil of rays is summed up by the contribution of both absorption and scattering. It is also known as ``extinction''.
\end{itemize}

The equation, $\sigma_{t} = \sigma_{a}+\sigma_{s}$, describes how the overall attenuation  $\sigma_{t}$  is influenced by both absorption $\sigma_{a}$ and scattering $\sigma_{s}$. In simpler terms, it breaks down the total impact on the medium into the portions caused by absorption and scattering.

The functions $f$ and  $u_{\uin}$  are the source and external radiative intensity. 
The  integral operator $ \mathcal{K}$ in (\ref{SRte1}) is given by
\[
\mathcal{K}u(\x,\s)=\int_{\Sc}\Phi(\x,\s, \s')u(\x,\s')d\sigma (\s'),\quad \int_{\Sc}\Phi(\x,\s, \s')d\sigma (\s') \equiv 1,\,\,  \forall \, \x \in \Omega,\,\,\s \in \Sc,
\]
where $\Phi$ is a non-negative normalized phase function, which describes the scattering property of the participating medium. 
The precise form of the phase function
is usually unknown for applications.
In this study, we consider a benchmark choice the Henyey-Greenstein (H-G) phase function which is  given by 
\begin{equation}\label{hgfun}
\Phi(\s, \s')=\dfrac{ 1}{2^{d-1} \pi}\dfrac{1-\eta^{2}}{\left(1+\eta^{2}-2\eta (\s \cdot \s')\right)^{3/2}},\,\,d=2,3.
\end{equation}
where the anisotropy factor $\eta \in (-1,1)$.

We assume that data of the model problem (\ref{SRte1}) satisfies the following conditions:
\begin{equation}\label{datarte}
\left\{
\begin{array}{ll}
\sigma_{t},\,\sigma_{s} \in L^{\infty}(\Omega),\quad  \sigma_{t}-\sigma_{s} \geq \sigma_{0}>0,\,\,\mbox{in}\,\,\Omega, \\
f(\x,\s) \in L^{2}(\Omega \times \Sc), \,\, f(\x, \cdot) \in \mathcal{C}(\Sc),\\
u_{\uin} \in L^{2}(\Gamma_{-}),
\end{array}
\right.
\end{equation}
where $\mathcal{C}(\Sc)$ is a continuous space over $\Sc$. Next, we denote the space $H^{1}_{2}(\Omega \times \Sc)$ by
\[
H^{1}_{2}(\Omega \times \Sc)=\{v \in L^{2}(\Omega \times \Sc),\,\s \cdot \nabla v \in L^{2}(\Omega \times \Sc)\}.
\]
with the norm 
\[
\|v\|_{H^{1}_{2}(\Omega \times \Sc)} = \left(\|v\|^{2}_{L^{2}(\Omega \times \Sc)}+\|\s \cdot \nabla v\|^{2}_{L^{2}(\Omega \times \Sc)}\right)^{1/2}.
\]
By using  (\ref{datarte}), it can be shown that the model problem has a unique solution $u \in H^{1}_{2}(\Omega \times \Sc)$. For more details, one may see \cite[Theorem 1]{sheng2021spherical}.

\bigskip
It is challenging to address the numerical approximation of the model problems because of the higher-dimensional RTE. In literature, the numerical methods for RTE can be categorized into the stochastic and deterministic approaches. In the stochastic approach, one of the major computational techniques is the  Monte Carlo simulations. There are various articles, that implemented the   Monte-Carlo simulation, for example, \cite{gentile2001implicit,howell,kong2008efficient,mcclarren2009modified}. Because of its iterative nature, the Monte-Carlo approach is not very time-efficient for the numerical approximation of RTE.  In the deterministic approach, some popular numerical methods are the Spherical harmonics ($P_N$)-method and the discrete-ordinate method (DOM), also called the $S_N$-method. The spherical harmonics method is based on the expansion of radiation intensity in angles using spherical harmonics basis functions. Detailed analysis of these methods for RTE can be found in 
\cite{egger2012mixed,Egger_SIAM16,egger_SIAM1,sheng2021spherical} and references therein. In DOM, the basic idea is to use a finite number of discrete directions to discretize the angular variables of RTE, which is also our area of interest.   DOM reduces the original RTE into a semi-discrete first-order hyperbolic system. The detailed description of this numerical approach can be found in \cite{koch1995discrete,lewis1984computational,kanschatnumerical}.

Owing to RTE's hyperbolic nature, one of the most widely used numerical methods for spatial discretization in the context of finite element schemes is the discontinuous Galerkin (DG) approach. DG methods for differential mathematical models have been extensively studied in the last few decades, see e.g. \cite{arnold2002unified,brezzi2004discontinuous,houston2002discontinuous} and references therein.  Recently, the local discontinuous Galerkin (LDG) finite element method for the fractional Allen-Cahn equation is studied in \cite{wang2023allen}. In \cite{yang2024optimal}, for the numerical simulation of the nonlinear fractional hyperbolic wave model, a second-order time discrete algorithm with a fast time two-mesh mixed finite element scheme is discussed.
DG methods for spatial discretization of RTE based on the DOM method are discussed in \cite{han} and also error analysis is presented. Later, a stabilized DG method along with DOM  is used for numerical simulation in   \cite{wang}. A more sophisticated numerical method comprises a multigrid preconditioner for $S_N$-DG approximation is presented in \cite{kanschat_SIAMJSC14}.
Recently, the upwind DG method for scaled RTE is discussed in \cite{sheng2021uniform}.  In \cite{peng2020low}, a low-rank numerical method for time-dependent RTE is established. A sparse grid DOM discontinuous Galerkin method for the RTE is investigated in \cite{Jianguo2021}. In \cite{ganesan2021operator}, we have proposed an operator-splitting based finite element method for the time-dependent RTE, where the key idea is to split the RTE model problem concerning the internal (angular) and the external (spatial) directions, resulting in a transient transport problem and a time-dependent integro-differential equation. A numerical simulation based on A reduced basis method for the stationary RTE model is proposed in \cite{peng2022reduced}.

\vskip 2mm

Discrete-ordinate DG type finite element methods for RTE are developed in \cite{han, wang}. In this article, we discuss the weak Galerkin (WG) method for the spatial discretization of the resulting semi-discrete system. One of the highlights of the WG method is that it has the advantage of the general mesh, which provides a high degree of versatility in numerical simulation, especially in the complex geometry domain. Unlike the DG method, the WG method is parameter-free while it has the same degrees of freedom as DG possesses.  More details of this method can be found in \cite{lin2018weak,mu2015weak,wang2014weak} and references therein. 
In the past, the WG finite element method has been used to solve several models. For example, the Cahn-Hilliard model is solved in \cite{wang2019weak}, and Biot's consolidation model is investigated in \cite{hu2018weak}.  Recently, an explicit WG method has been used for numerical simulation of first-order hyperbolic systems in \cite{zhang2022explicit}.
Inspired by these works, we explore the possibility of using the WG method for the RTE model, where we will use a modified version of the weak Galerkin method introduced in  \cite{wang2014modified,gao2020modified} for spatial discretization. 

\vskip 2mm 
This is the first article that focuses on discretizing the spatial direction of RTE utilizing the WG methods. The proposed numerical method combines DOM for the angular variable and the WG method for the spatial variable. The resulting technique is called the discrete-ordinate weak Galerkin method or DOWG method for short. The stability result of the DOWG method is derived. A \emph{priori} error estimate is obtained by combining the projection and discretization errors. 

\vskip 2mm
The rest of the research is framed as follows: In Section \ref{Sec2}, we define the DOM for the model problem (\ref{SRte1}). Then, we introduce the WG method for the resulting semi-discrete hyperbolic system and this section ends with the stability result of the proposed  DOWG method. Section \ref{Sec3} discusses the error estimate of the DOWG method by combining the error estimate associated with both the numerical methods,  \emph{i.e}, DOM, and WG method. The numerical experiments that support our theoretical findings are covered in Section \ref{Sec4} and the article is summarized in Section \ref{Sec5}.

\vskip 2mm
\noindent
{\it Notations:}  We have used the standard notations for Sobolev space and norm/semi-norm 
\cite{brnnerfem} throughout this article. And 
$C$, sometimes with a subscript,   denotes a generic positive constant,
the value of which can be different in different occurrences.

\section{Discrete-ordinate WG finite element method}\label{Sec2}
In this section, to solve the RTE, a discrete-ordinate weak Galerkin finite element method is presented (\ref{SRte1}). The numerical simulations are developed in two phases: First, the
discrete-ordinate approach is applied to discretize the angular variable, which deduces the model problem (\ref{SRte1})
into a system of hyperbolic PDE.
Then,  the resulting system of differential equations is further discretized by the weak Galerkin method.

\subsection{Discrete-ordinate method}

Before discussing the discrete-ordinate method, we briefly go through the quadrature rule, which is instrumental for the proposal of the discrete method. For any continuous function $F$, we can write due to the quadrature approximation
\begin{equation}\label{gen_quad}
    \int_{\Sc}F(\s)d\sigma(\s) \approx \sum_{m=0}^{M}\w_{m}F(\s_m), \quad \w_m > 0,\,\,\ \s_m \in \Sc,\, 0\leq m \leq M.
\end{equation}

For $\Sc \in \mathbb{R}^{2}$, the quadrature rule (\ref{gen_quad}) is expressed as 
\begin{equation}\label{gen_quad_2d}
    \int_{\Sc}F(\s)d\sigma(\s) = \int_{ 0}^{2\pi}\overline{F}(\theta)d\theta \approx \sum_{m=0}^{M}\w_{m}\overline{F}(\theta_m), \quad \w_m > 0,\,\,\ \theta_m \in [0, 2\pi],\, 1\leq m \leq M,
\end{equation}
where the angular variable is written in the spherical coordinate form. Afterwards, the standard quadrature formula is used. More details can be found in any standard numerical analysis book, for example, see \cite{suli2003introduction}. In this paper, we have applied the composite trapezoidal quadrature formula in (\ref{gen_quad_2d}). To achieve this, we adopted the angular spacing parameter  $h_{\theta}=2\pi/M$ with weights $\w_{0}=\w_{M}=\dfrac{h_{\theta}}{2}$ with $\w_{m} = h_\theta$ and the  discrete points $\theta_m = mh_{\theta}$, for all $1 \leq m \leq M-1$.

Similarly, for $\Sc \in \mathbb{R}^{3}$, the following quadrature rule is used in (\ref{gen_quad}) to obtain  
\begin{equation}\label{gen_quad_3d}
\begin{array}{ll}
\displaystyle  \int_{\Sc}F(\s)d\sigma(\s) &
\displaystyle = \int_{ 0}^{2\pi} \int_{ 0}^{\pi}\overline{F}(\theta, \psi) \sin\theta d\theta d\psi, \quad (\theta, \psi) \in [0, \pi] \times [0,2\pi] \\[12pt] 
  & \displaystyle \approx \dfrac{\pi}{m} \sum_{j=1}^{2m}\sum_{i=1}^{m}\w_{i}\overline{F}(\theta_i, \psi_j), \quad \w_i > 0,
\end{array}
\end{equation}
where $\{\theta_i\}$  are chosen so that $\{\cos \theta_i\}$ and $\{\w_i\}$ are the Gauss-Legendre nodes and weights on $[-1,1]$. The discrete points $\{\psi_j\}$ are evenly spaced on $[0,2\pi]$ with a spacing of $\pi/m$. 
By following the quadrature error estimates  \cite{hesse2006cubature},   we have 
\begin{equation}\label{quad3d_err}
    \left|\int_{\Sc}F(\s)d\sigma(\s) - \sum_{m=0}^{M}\w_{m}F(\s_m)\right| \leq C_{r} n^{-r} \|F\|_{r, \Sc},\,\, F \in H^{r}(\Sc \in \mathbb{R}^{3}), \,\, r \geq 1,
\end{equation}
where $C_r$ is a positive constant depending only on $r$,  and $n$ denotes the degree of precision of the quadrature. Here,  norm  $\|\cdot\|_{r, \Sc}$ in (\ref{quad3d_err})
\begin{equation}\label{hnorm}
\|v\|^{2}_{r, \Sc}:=\int_{\Sc}\sum_{k =0}^{r}\left|\nabla^{(k)}_{\s} v\right|^{2}, \qquad \nabla^{(k)}_{\s} v := \left\{\nabla^{(\alpha)}_{\s} v \big| |\alpha| = k\right\},
\end{equation}
where  $\alpha$ is a multi-index and $|\nabla^{(k)}_{\s} v| =  \sum_{|\alpha|=k}\left(|\nabla^{(\alpha)}_{\s} v|^{2}\right)^{1/2}$. 
By using the normalized property of the scattering phase function with the estimate (\ref{quad3d_err}), for any fixed $\x \in \Omega$, we get
\begin{equation}\label{quad01}
\left|1- \sum_{\ell=0}^{M}\w_{m}\Phi(\s_{m}\cdot \s_{\ell})\right| \leq C_{r} n^{-r}\|\Phi(\s_{m}\cdot)\|_{r,
	\Sc},\quad r \geq 1.
\end{equation}

Using the above numerical quadrature formulas, the integral operator  $ \mathcal{K} $ is approximated  in the following unified way:
\begin{equation}\label{SRte2}
\mathcal{K}u(\x,\s) \approx \mathcal{K}_{d} u(\x, \s_{m}) = \sum_{m=0}^{M}\w_{m}\Phi(\x,\s\cdot \s_{m})u(\x, \s_{m}),\quad \w_{m} >0, \,\,\mbox{and}\,\,\s_{m} \in \Sc,\,\,0 \leq  m \leq M. 
\end{equation}

\subsubsection{Semi-discrete problem}
By employing the approximation (\ref{SRte2}), the continuous problem (\ref{SRte1}) is  discretized in each angular direction $\s_m$ as follows:
\begin{equation}\label{SRtde1}
\left\{
\begin{array}{ll}
\displaystyle \s_{m} \cdot \nabla u^{m}(\x)+\sigma_{t}(\x) u^{m}(\x)-\sigma_{s} \sum_{\ell=0}^{M}\w_{\ell}\Phi(\x,\s_{m}\cdot \s_{\ell})u^{\ell}(\x) =f^{m}(\x),\quad  \mbox{in} \,\, \Omega,\\[12pt]
u^{m}(\x)=u_{\uin}^{m}(\x),\quad  \mbox{on} \,\, \partial \Omega_{m}^{-} \equiv \partial \Omega_{\s_{m}}^{-}, \qquad \qquad  \qquad  0 \leq m \leq M,
\end{array}
\right.
\end{equation}
where $f^{m}(\x)=f(\x,\s_{m}),\,u_{\uin}^{m}=u_{\uin}(\x,\s_{m})$ and the spatial boundaries $\Omega_{\s_m}^{\pm}$ are
\[
\partial \Omega_{\s_m}^{-}=\{\x \in \partial \Omega: \s_m \cdot \mathbf{n}(\x)<0\},\quad \partial \Omega_{\s_m}^{+}=\{\x \in \partial \Omega: \s_m \cdot \mathbf{n}(\x) \geq 0\},
\]
and the incoming and outgoing boundaries are given by
$\Gamma^{\pm}_{m}=\{(\x,\s_m):\x \in \partial \Omega_{\s}^{\pm}\}$.

Note that, the  exact solution $u^{m}=u^{m}(\x)$ of the semi-discrete  system (\ref{SRtde1}) is an approximation of the model solution $u$ of the continuous problem (\ref{SRte1}) at points $(\x, \s_{m})$ for $m=0,2,\ldots,M$.

Before going into the details of the error estimate, consider the following regularity conditions for the solution of the discrete method (\ref{SRtde1}):
In the $2D$ setting, we assume that the exact solution $u$ of the continuous problem (\ref{SRte1}) and the phase function satisfy
\begin{equation}\label{regasm2d}
\left\{
\begin{array}{ll}
u^{m} \in H^{p}(\Omega),\quad \mbox{for }\,\, p \geq 0,\quad 1 \leq m \leq M,\\[6pt]
\displaystyle u \in L^{2}(\Omega, \mathcal{C}^{2}(\Sc)) \qquad \mbox{and} \quad\sup_{(\x,\s) \in \Omega \times \Sc}\|\Phi''(\x,\s \cdot)\|_{0,\infty,\Sc} < C_{\Phi}, \quad \mbox{for }\,\, r \geq 1, \quad \mbox{where}\,\, \Phi''(\x,t) = \dfrac{d^{2} \Phi}{d t^{2}}.
\end{array}
\right.
\end{equation}
And, in the $3D$ setting
\begin{equation}\label{regasm}
\left\{
\begin{array}{ll}
u^{m} \in H^{p}(\Omega),\quad \mbox{for }\,\, p \geq 1,\quad 1 \leq m \leq M,\\[6pt]
\displaystyle u \in L^{2}(\Omega,H^{r}(\Sc)) \qquad \mbox{and} \quad \sup_{(\x,\s) \in \Omega \times \Sc}\|\Phi(\x,\s \cdot)\|_{r,\Sc} < C_{\Phi_*}, \quad \mbox{for }\,\, r \geq 1.
\end{array}
\right.
\end{equation}
where the Bochner spaces  $L^{2}(\Omega, \mathcal{C}^{2}(\Sc))$ and $L^{2}(\Omega,H^{r}(\Sc))$ are given  by
\begin{align*}
 L^{2}(\Omega,\mathcal{C}^{2}(\Sc))=\left\{v :\Omega \to \mathcal{C}^{2}(\Sc)): \int_{\Omega}\|v\|_{\mathcal{C}^{2}(\Sc))}< \infty \right\},\\[8pt]  
 L^{2}(\Omega,H^{r}(\Sc)=\left\{v :\Omega \to H^{r}(\Sc): \int_{\Omega}\|v\|_{r,\Sc}< \infty \right\}.
\end{align*}

For simplification of the presentation, we introduce
\[
b(\x)=\max_{1 \leq \ell \leq M} \sum_{\ell=0}^{M}\w_{m}\Phi(\x,\s_{\ell}\cdot \s_{m}).
\]
Then, from  (\ref{datarte}) and (\ref{quad01}), one can deduce that
\begin{equation}\label{quadasum}
\sigma_{t} (\x) -\sigma_{s}(\x)  b(\x) \equiv  (\sigma_{a}+\sigma_{s}) (\x) -\sigma_{s}(\x)  b(\x)  \geq C_{*} \quad \mbox{in} \,\,\Omega,\quad C_{*} \,\,\mbox{a positive constant}.
\end{equation}
This condition will be utilized while deriving the stability and convergence result of the numerical approximation.

In the next few lines, we discuss the error estimate using the discrete-ordinate method (\ref{SRtde1}).  
\begin{theorem}\label{domerrthm}
Let  $u$  and  $u^m$ are the solutions of the model problem (\ref{SRte1}) and the semi-discrete problem (\ref{SRtde1}), respectively.
In the $2D$ case,   
 \begin{equation}\label{domerror}
\left(\sum_{m=0}^{M}\w_{\ell}\sum_{T \in \mathcal{T}_{h}}\|u^{m}(\x)-u(\x,\s_{m})\|^{2}_{T} \right)^{1/2} = \mathcal{O}(h_{\theta}^2),
\end{equation}
where $h_{\theta}$ is sufficiently small. And, in the $3D$ case, we have
 \begin{equation}\label{domerror1}
\left(\sum_{m=0}^{M}\w_{\ell}\sum_{T \in \mathcal{T}_{h}}\|u^{m}(\x)-u(\x,\s_{m})\|^{2}_{T} \right)^{1/2}\leq C\,n^{-r}\left(\int_{\Omega}\|u\|^{2}_{r,\Sc}\,\mbox{d}\x\right)^{1/2},
\end{equation}
where  $C$ is positive constant depending on $r$ and the phase function $\Phi$.
\end{theorem}
{\bf Proof.} By using the regularity properties (\ref{regasm2d})-(\ref{regasm}), the proof of this theorem is followed from \cite[Theorem 4.3, Theorem 4.4]{wang}. \hfill \eop

\subsection{Weak Galerkin finite element method}
Let $\mathcal{T}_{h}$ be a mesh partition of the spatial  domain $\Omega$ consisting of elements that are closed and simply connected polygons in $\mathbb{R}^{2}$ (polyhedrons in $\mathbb{R}^{3}$)
satisfying a set of shape-regular conditions \cite{mu2015weak}. Let $\mathcal{E}_{h}$ is the set of all edges/faces in $\mathcal{T}_{h}$ with $\mathcal{E}^{0}_{h}=\mathcal{E}_{h} \setminus \partial \Omega_{h}$ is the set of all interior boundaries (faces for $d = 3$ or edges for $d = 2$). And  $\partial \Omega_{h}$ be the subset of  $\mathcal{E}_{h}$ of all boundaries on $\partial \Omega$. For every element $T \in \mathcal{T}_{h}$, we denote by $h_T$ the diameter
of $T$ and the mesh size $h= \max_{T \in \mathcal{T}_{h}}
h_{T}$ for $\mathcal{T}_{h}$.

For a fixed angular direction $\s_{m}$,  the incoming and outgoing boundaries of the mesh element $T \in \mathcal{T}_{h}$ is given by
\[
\partial^{-}T=\{\x \in \partial T: \s_{m} \cdot \mathbf{n}(\x)<0\}, \qquad \partial^{+}T=\{\x \in \partial T: \s_{m} \cdot \mathbf{n}(\x) \geq 0\}.
\]

Let $T_{1}$ and $T_{2}$ are the mesh elements sharing a common edge $ e $. We define the average of the scalar valued  $ v $ and  vector valued functions $\mathbf{v}$ on edge $ e $ as follows:
\[
\{v\}_{e}=
\left\{
\begin{array}{ll}
\dfrac{1}{2}(v|_{T_1}+v|_{T_2}), & e \in \mathcal{E}^{0}_{h},\\[8pt]
v, & e \in \partial \Omega,
\end{array}
\right.
\quad 
\{\mathbf{v}\}_{e}=
\left\{
\begin{array}{ll}
\dfrac{1}{2}(\mathbf{v}|_{T_1}+\mathbf{v}|_{T_2}), & e \in \mathcal{E}^{0}_{h},\\[8pt]
\mathbf{v}, & e \in \partial \Omega.
\end{array}
\right.
\]
And, deﬁne the jump of $v$ on $e$ by $[v] = v^{+}-v^{-}$.

For convenience, we introduce the following notations:
\[
\displaystyle (\phi,\psi)_{\mathcal{T}_{h}}=\sum_{T \in \mathcal{T}_{h}}\int_{T} \phi \psi \mbox{d}\x, \,\, \big\langle \phi,\psi \rangle_{\partial\mathcal{T}_{h}}=\sum_{T \in \mathcal{T}_{h}}\int_{\partial T} \phi \psi\, \mbox{ds},\,\, \big\langle \phi,\psi \rangle_{\partial \Omega_{h}}=\sum_{T \in \mathcal{T}_{h}}\int_{\partial T \cap \partial \Omega_{h}} \phi \psi\, \mbox{ds}. 
\]

Now, we introduce the finite element space $V_{h}$.
Given an integer $k \geq 1$,  the finite element space is defined as follows
\[
V_{h}=\{v_{h} \in L^{2}(\Omega): v_{h}|_{T} \in  \mathbb{P}_{k}(T),\,\,\forall \,T \in \mathcal{T}_{h}\},
\]
where $ \mathbb{P}_{k}(T) $ is the space of polynomials on set $T$ with degree no more than $k$.

\begin{definition}
For any function $v^{m}_{h}  \in V_h$, on each element $T$, we define the weak gradient operator $\nabla_{w} v_{h}^{m}|_{T} \in [\mathbb{P}_{k-1}(T)]^{d}$ as follows:
\begin{equation}\label{def:weakgrad}
 (\nabla_{w} v_{h}^{m},\mathbf{q}_{h})_{T}=-(v^{m}_{h},\nabla \cdot \mathbf{q}_{h})_{T}+\big\langle  \{v^{m}_{h}\},\mathbf{q}_{h} \cdot \mathbf{n} \big\rangle_{\partial T},\quad \forall \mathbf{q}_{h} \in [\mathbb{P}_{k-1}(T)]^{d}.   
\end{equation}
And we define the weak divergence $\s_{m} \cdot \nabla_{w} v^{m}_{h} \in \mathbb{P}_{k}(T)$ associated with $\s_{m} \cdot \nabla v_{h}^{m}$ on $T$, as the unique polynomial satisfying 
\begin{equation}\label{def:weakdiv}
(\s_{m} \cdot \nabla_{w} v^{m}_{h},w_{h})=-(v^{m}_{h},\s_{m} \cdot \nabla w_{h})+\big \langle  \{v^{m}_{h}\},\s_{m} \cdot \mathbf{n}\, w_{h} \big \rangle_{\partial T},\quad \forall w_{h} \in \mathbb{P}_{k}(T).
\end{equation}
\end{definition}

Next,  we introduce the projection operators, which will be used in the error analysis of the weak Galerkin method. For each element $ T \in \mathcal{T}_{h}$, let $\mathcal{Q}_{T}:L^{2}(T) \to \mathbb{P}_{k}(T)$
and $ \Pi_{T} :[L^{2}(T)]^{d} \to [\mathbb{P}_{k-1}(T)]^{d}$
be the $L^2$-projections on the associated local polynomial spaces respectively. For any mesh element $T$, we define  the global projections $ \mathcal{Q}_{h} $ and $ \Pi_{h} $ element-wise by
\[
\begin{array}{ll}
(\mathcal{Q}_{h}v)|_{T}= \mathcal{Q}_{T}(v|_{T}),& \forall \, v \in L^{2}(T),\\[4pt]
(\Pi_{h}\mathbf{q})|_{T}=\Pi_{T}(\mathbf{q}|_{T}),& \forall \, \mathbf{q} \in [L^{2}(T)]^{d}.
\end{array}
\]

\vskip 2mm
The finite dimensional space $\mathbf{V}_{h}$ is  defined as   $\mathbf{V}_{h}=(V_{h})^{M}$ for the semi-discrete problem (\ref{SRtde1}) and   an  element  $\mathbf{v}_{h} \in \mathbf{V}_{h} $ is given by $\mathbf{v}_{h}:=\{v_{h}^{m}\}_{m=0}^{M}$.

\vskip 2mm
The  WG finite element method   for the semi-discrete problem (\ref{SRtde1}) reads:
Find $u_{h}^{m} \in V_{h}$  such that
\begin{equation}\label{SRtfe1}
\displaystyle \A_{c}(u_{h}^{m},v_{h}^{m})+\A_{st}(u_{h}^{m},v_{h}^{m})=\F^{\s}(v_{h}^{m}), \quad \forall \,v_{h}^{m} \in V_{h},
\end{equation}
where $\A_{c}(\cdot,\cdot),\,\A_{st}(\cdot,\cdot)$ and $\F^{\s}(\cdot)$ are given by
\[
\begin{array}{ll}
\displaystyle \A_{c}(u_{h}^{m},v_{h}^{m})=\displaystyle \sum_{T \in \mathcal{T}_{h}}\bigg(\s_{m} \cdot \nabla_{w} u_{h}^{m}+\sigma_{t} u_{h}^{m}-\sigma_{s} \sum_{\ell=0}^{M}\w_{\ell}\Phi(\x,\s_{m}\cdot \s_{\ell})u_{h}^{\ell},v_{h}^{m}\bigg)_{T}-\langle \s_{m} \cdot \mathbf{n} \,u_{h}^{m},v_{h}^{m}\rangle_{\partial \Omega^{-}_{h}},\\[16pt]
\A_{st}(u_{h}^{m},v_{h}^{m})=\displaystyle \sum_{T \in \mathcal{T}_{h}}\left\langle \s_{m} \cdot \mathbf{n} \big(u^{m}_{h}-\{u^{m}_{h}\}\big),v^{m}_{h}-\{v^{m}_{h}\}\right\rangle_{\partial^{+} T},\\[16pt]
\F^{\s}(v_{h}^{m})=(f^{m}, v_{h}^{m})_{\mathcal{T}_{h}}-\langle \s_{m} \cdot \mathbf{n} \,u_{\uin}^{m},v_{h}^{m}\rangle_{\partial \Omega^{-}_{h}},
\end{array}
\]
where, the incoming boundary
$\partial \Omega^{-}_{h} = \left\{ e  \in \partial \Omega_{h}  : \s \cdot \mathbf{n}(\x)<0,\,\x \in e\right\}$. Note that, we have imposed the boundary condition in the weak sense. 

\subsection{Fully-discrete method}
The global formulation of the DOWG method (\ref{SRtfe1}) is  expressed as:
Find $\bu_{h} =\{u_{h}^{m} \} \in \mathbf{V}_{h}$ such that
\begin{equation*}
\begin{array}{ll}
\displaystyle \sum_{m=0}^{M}\w_{m} \, \A_{c}(u_{h}^{m},v^{m}_{h})+\sum_{m=0}^{M}\w_{m} \,\A_{st}(u_{h}^{m},v^{m}_{h})=\displaystyle \sum_{m=0}^{M}\w_{m} \F^{\s}(v_{h}^{m}).
\end{array}
\end{equation*}

In the simplified version, the  DOWG method  is presented as follows:
\begin{equation}\label{wg1}
\left\{
\begin{array}{ll}
\mbox{Find}\,\, \bu_{h} \in \mathbf{V}_{h} \,\,\mbox{such that}\\[6pt]
\A_{WG}(\bu_{h},\bv_{h})=\F(\bv_{h}), \,\,\forall \, \bv_{h} \in \mathbf{V}_{h},
\end{array}
\right.
\end{equation}
where the bilinear form $\A_{WG}(\cdot,\cdot)$ and the linear functional $\F^{\s}(\cdot)$ are given by
\[
\A_{WG}(\bu_{h},\bv_{h})=\sum_{m=0}^{M}\w_{m} \big(\A_{c}(u_{h}^{m},v^{m}_{h})+\A_{st}(u_{h}^{m},v_{h}^{m})\big),\quad \F(\bv_{h})=\sum_{m=0}^{M}\w_{m} \F^{\s}(v_{h}^{m}).
\]
For any $\bv_{h} \in \mathbf{V}_{h}$ , the proposed  norm $\vertiii{\cdot}$ is  defined as
\begin{equation}\label{normdef}
    \vertiii{\bv_{h}}^{2}= \sum_{m=0}^{M}\w_{m}\vertiii{v^{m}_{h}}_{\mathcal{T}_{h}}^{2}, 
\end{equation}
where  $ \vertiii{\cdot}_{\mathcal{T}_{h}} $ is given by
\[
 \vertiii{v_{h}^{m}}_{\mathcal{T}_{h}}^{2}=\sum_{T \in \mathcal{T}_{h}} \left( \|v^{m}_{h}\|^{2}_{T}+\left\| |\s_{m} \cdot \mathbf{n}|^{1/2}\big(v_{h}^{m}-\{v_{h}^{m}\}\big)\right\|^{2}_{\partial T}\right)+\| |\s_{m} \cdot \mathbf{n}|^{1/2}v_{h}^{m}\|^{2}_{\partial \Omega_{h}}. 
\]

\section{Error analysis}\label{Sec3}
In this section, we discuss the error estimate of the DOWG method. Firstly, we start with the stability result of the proposed numerical method (\ref{wg1}). By using this stability result, we derive the error estimate of the numerical solution. 
\subsection{Stability result}
Now, we discuss the coercivity of the bilinear form $ \A_{WG}(\cdot,\cdot) $ in the following lemma.
\begin{lemma}\label{stab}
	For sufficiently small $h$, we have
	\[
\alpha	\vertiii{\bv_{h}}^{2} \leq  \, \A_{WG}(\bv_{h},\bv_{h}), \,\,\forall \, \bv_{h} \in \mathbf{V}_{h},
	\]
	where $\alpha$ is a  positive constant.
\end{lemma}
{\bf Proof.} 
By using the definition of  weak divergence (\ref{def:weakdiv}) and integral by parts, it can be easily deduced that
\[
\left(\s_{m}\cdot \nabla_{w}v^{m}_{h}, v^{m}_{h}\right)_{\mathcal{T}_{h}}=-\dfrac{1}{2}\left\langle \s_{m} \cdot \mathbf{n}\big(v_{h}^{m}-\{v_{h}^{m}\}\big),v_{h}^{m}-\{v_{h}^{m}\}\right\rangle_{\partial \mathcal{T}_{h}}+\dfrac{1}{2}\langle \s_{m} \cdot \mathbf{n} \,v_{h}^{m},v_{h}^{m}\rangle_{\partial \Omega_{h}}.
\]
Then, we have 
\[
\begin{array}{ll}
\A_{c}(v_{h}^{m},v_{h}^{m})+\A_{st}(v_{h}^{m},v_{h}^{m})&=\displaystyle\sum_{T \in \mathcal{T}_{h}}\left(\sigma_{t} v^{m}_{h} , v^{m}_{h}\right)_{T}+I+\dfrac{1}{2}\langle \s_{m} \cdot \mathbf{n} \,v_{h}^{m},v_{h}^{m}\rangle_{\partial \Omega_{h}}\\[16pt]
&\quad -\dfrac{1}{2}\left\langle \s_{m} \cdot \mathbf{n}\big(v_{h}^{m}-\{v_{h}^{m}\}\big),v_{h}^{m}-\{v_{h}^{m}\}\right\rangle_{\partial \mathcal{T}_{h}}\\[16pt]
&\displaystyle \quad +\sum_{T \in \mathcal{T}_{h}}\left\langle \s_{m} \cdot \mathbf{n} \big(v^{m}_{h}-\{v^{m}_{h}\}\big),v^{m}_{h}-\{v^{m}_{h}\}\right\rangle_{\partial^{+} T},
\end{array}
\]
where $I$ is given by
\[
I=-\sum_{T \in \mathcal{T}_{h}}\left(\sigma_{s} \sum_{\ell =0}^{M}\w_{\ell}\Phi(\x,\s_{m}\cdot \s_{\ell})v^{\ell}_{h} , v^{m}_{h}\right)_{T},
\]
and it  satisfies
\[
|I| \leq \sum_{T \in \mathcal{T}_{h}} \left(\sigma_{s} bv^{m}_{h} , v^{m}_{h}\right)_{T}.
\] 
By using condition (\ref{quadasum}), we obtain that
\[
\begin{array}{ll}
\A_{c}(v_{h}^{m},v_{h}^{m})+\A_{st}(v_{h}^{m},v_{h}^{m})&\geq \displaystyle  \left((\sigma_{t}-\sigma_{s}b) v^{m}_{h} , v^{m}_{h}\right)+\dfrac{1}{2}\langle \s_{m} \cdot \mathbf{n} \,v_{h}^{m},v_{h}^{m}\rangle_{\partial \Omega_{h}}\\[14pt]
& \quad +\dfrac{1}{2}\left\langle \s_{m} \cdot \mathbf{n}\big(v_{h}^{m}-\{v_{h}^{m}\}\big),v_{h}^{m}-\{v_{h}^{m}\}\right\rangle_{\partial \mathcal{T}_{h}}\\[14pt]
&\geq \left(C_{*} v^{m}_{h} , v^{m}_{h}\right)+\dfrac{1}{2}\langle \s_{m} \cdot \mathbf{n} \,v_{h}^{m},v_{h}^{m}\rangle_{\partial \Omega_{h}}\\[14pt]
& \quad +\dfrac{1}{2}\left\langle \s_{m} \cdot \mathbf{n}\big(v_{h}^{m}-\{v_{h}^{m}\}\big),v_{h}^{m}-\{v_{h}^{m}\}\right\rangle_{\partial \mathcal{T}_{h}}.
\end{array}
\] 
Therefore, the bilinear form satisfies
\[
\begin{array}{ll}
\A_{WG}(\bv_{h},\bv_{h}) & \displaystyle \geq \alpha \sum_{m=0}^{M}\w_{m}\vertiii{v^{m}_{h}}_{\mathcal{T}_{h}}^{2}\\ [12pt]
& = \alpha	\vertiii{\bv_{h}}^{2},\qquad  \alpha=\min\{C_{*},\frac{1}{2}\}.
\end{array}
\]
This completes the proof.
\hfill
\eop 

\vskip 2mm
The direct consequence of the above lemma is the uniqueness result.
\begin{theorem}\label{Uniqe}
There exists a unique solution of the DOWG method (\ref{wg1}).
\end{theorem}

\vskip 2mm
Next, we state the trace inequality: For any function $\phi \in H^{1}(T)$, the following trace inequality holds
true (see \cite{wang2014weak}):
\begin{equation}\label{traceeq}
\|\phi\|^{2}_{\partial T} \leq C \left(h_{T}^{-1}\|\phi\|^{2}_{ T}+h_{T}\|\nabla \phi\|^{2}_{ T}\right).
\end{equation}

\subsection{Error estimate}
Next, we discuss the error estimate of the DOWG method. Firstly, we establish an error equation associated with the discrete problem (\ref{wg1}). 

Define the components $\mathcal{L}_{j}(\cdot,\cdot)$ as follows:
\begin{align}\label{L_comp}
    \mathcal{L}_{1}(u^{m},v_{h}^{m})&=(u^{m}-\mathcal{Q}_{h}u^{m},\s_{m}\cdot \nabla v_{h}^{m})_{\mathcal{T}_{h}},\nonumber\\[8pt]
 \mathcal{L}_{2}(u^{m},v_{h}^{m})&=-\big\langle u^{m}-\{\mathcal{Q}_{h}u^{m}\},\s_{m} \cdot \mathbf{n} (v_{h}^{m}-\{v_{h}^{m}\})\big \rangle_{\partial \mathcal{T}_{h}},\\[8pt]
 \mathcal{L}_{3}(u^{m},v_{h}^{m})&= \displaystyle -\left(\sigma_{t} (u^{m}-\mathcal{Q}_{h}u^{m})-\sigma_{s} \sum_{\ell =0}^{M}\w_{\ell}\Phi(\x,\s_{m}\cdot \s_{\ell})(u^{\ell}-\mathcal{Q}_{h}u^{\ell}) , v_{h}^{m}\right)_{\mathcal{T}_{h}},\nonumber\\[8pt]
 \mathcal{L}_{4}(u^{m},v_{h}^{m})&=\big\langle \s_{m} \cdot \mathbf{n}(u^{m}-\mathcal{Q}_{h}u^{m}), v_{h}^{m}\big\rangle_{\partial \Omega^{-}_{h}}. \nonumber
\end{align}

\begin{lemma}\label{erreq1}
	The error term
	$\mathbf{e}_{\mathcal{Q}_{h}}=(\mathcal{Q}_{h}\bu-\bu_{h})$
	satisfies:
\begin{equation}\label{em0}
\A_{WG}(\mathbf{e}_{\mathcal{Q}_{h}},\bv_{h})=\mathbf{L}(\bu,\bv_{h}),
\end{equation}
with 
\[
\mathbf{L}(\bu,\bv_{h})=\sum_{\ell =0}^{M}\w_{m} \mathcal{L}(u^{m},v_{h}^{m}),\qquad \mathcal{L}(u^{m},v^{m}_{h})=\sum_{j=1}^{4}\mathcal{L}_{j}(u^{m},v_{h}^{m})+\A_{st}(\mathcal{Q}_{h}u^{m},v_{h}^{m}),
\]
where  $\mathcal{L}_{j}(\cdot,\cdot)$ are given in (\ref{L_comp}). 
\end{lemma}
{\bf Proof.} 
By using the integration by parts and definition (\ref{def:weakdiv}), we deduce that
\begin{equation}\label{er1}
\begin{array}{ll}
\left(\s_{m}\cdot \nabla u^{m},v^{m}_{h}\right)_{\mathcal{T}_{h}}&=-(u^{m},\s_{m}\cdot \nabla v_{h}^{m})_{\mathcal{T}_{h}}+\langle u^{m},\s_{m} \cdot \mathbf{n} v_{h}^{m} \rangle_{\partial \mathcal{T}_{h}}\\[8pt]
&=-(\mathcal{Q}_{h}u^{m},\s_{m}\cdot \nabla v_{h}^{m})_{\mathcal{T}_{h}}-(u^{m}-\mathcal{Q}_{h}u^{m},\s_{m}\cdot \nabla v_{h}^{m})_{\mathcal{T}_{h}}\\[8pt]
&
\quad +\big\langle \{\mathcal{Q}_{h}u^{m}\},\s_{m} \cdot \mathbf{n} v_{h}^{m} \big\rangle_{\partial \mathcal{T}_{h}}+\big\langle u^{m}-\{\mathcal{Q}_{h}u^{m}\},\s_{m} \cdot \mathbf{n} (v_{h}^{m}-\{v_{h}^{m}\}) \big\rangle_{\partial \mathcal{T}_{h}}
\\[8pt]
&=\left(\s_{m}\cdot \nabla_{w} (\mathcal{Q}_{h}u^{m}),v^{m}_{h}\right)_{\mathcal{T}_{h}}-\mathcal{L}_{1}(u^{m},v_{h}^{m})-\mathcal{L}_{2}(u^{m},v_{h}^{m}). 
\end{array}
\end{equation}
Next, one can  obtain that 
\begin{equation}\label{er2}
\left(\sigma_{t} u^{m}-\sigma_{s} \sum_{\ell =0}^{M}\w_{\ell}\Phi(\x,\s_{m}\cdot \s_{\ell})u^{\ell} , v^{m}_{h}\right)_{\mathcal{T}_{h}}=
\displaystyle \left(\sigma_{t} \mathcal{Q}_{h}u^{m}-\sigma_{s} \sum_{\ell =0}^{M}\w_{\ell}\Phi(\x,\s_{m}\cdot \s_{\ell})\mathcal{Q}_{h}u^{\ell} , v^{m}_{h}\right)_{\mathcal{T}_{h}} -\mathcal{L}_{3}(u^{m},v_{h}^{m})
\end{equation}
and
\begin{equation}\label{er02}
-\langle \s_{m} \cdot \mathbf{n} \,u^{m},v_{h}^{m}\rangle_{\partial \Omega^{-}_{h}}=
-\langle \s_{m} \cdot \mathbf{n} \,\mathcal{Q}_{h}u^{m},v_{h}^{m}\rangle_{\partial \Omega^{-}_{h}}-\mathcal{L}_{4}(u^{m},v_{h}^{m}). 
\end{equation}
Finally, by combining (\ref{er1})-(\ref{er02}), we get
\[
\A_{c}(\mathcal{Q}_{h}u^{m},v_{h}^{m})-\sum_{j=1}^{4}\mathcal{L}_{j}(u^{m},v_{h}^{m})=(f^{m}, v_{h}^{m})_{\mathcal{T}_{h}}-\langle \s_{m} \cdot \mathbf{n} \,u_{\uin}^{m},v_{h}^{m}\rangle_{\partial \Omega^{-}_{h}}.
\]
Adding $\A_{st}(\mathcal{Q}_{h}u^{m},v_{h}^{m})$ both side of the above equation to obtain 
\begin{equation}\label{er3}
\A_{c}(\mathcal{Q}_{h}u^{m},v_{h}^{m})+\A_{st}(\mathcal{Q}_{h}u^{m},v_{h}^{m})=(f^{m}, v_{h}^{m})_{\mathcal{T}_{h}}-\langle \s_{m} \cdot \mathbf{n} \,u_{\uin}^{m},v_{h}^{m}\rangle_{\partial \Omega^{-}_{h}}+\mathbf{L}(u^{m},v_{h}^{m}).
\end{equation}
The above equation into the global formulation is given by
\begin{equation}\label{er4}
\begin{array}{ll}
\A_{WG}(\mathcal{Q}_{h}\bu,\bv_{h})=\F(\bv_{h})+\mathbf{L}(\bu,\bv_{h}).
\end{array}
\end{equation}
By subtracting the discrete problem (\ref{wg1}) from (\ref{er4}), we obtain the required error equation.
\hfill \eop

In the next lemma,  we derive bounds for each term of $\mathbf{L}(\bu,\bv_{h})$. Then, we use the coercivity of the DOWG method to find the estimates of error term $\mathbf{e}_{\mathcal{Q}_{h}}$.
\begin{lemma}\label{errest0}
	 Let $ u^{m}$ be the solution of the discrete problem (\ref{SRtde1}). For $v^{m}_{h} \in V_{h}$, we have the following estimates:
\[
\begin{array}{ll}
|\mathcal{L}_{j}(u^{m},v_{h}^{m})|& \leq C h^{k+1}\|u^{m}\|_{k+1}\vertiii{v_{h}^{m}}_{\mathcal{T}_{h}},\,\,j=1,3,\\[8pt]
|\mathcal{L}_{j}(u^{m},v_{h}^{m})|& \leq C h^{k+1/2}\|u^{m}\|_{k+1}\vertiii{v_{h}^{m}}_{\mathcal{T}_{h}},\,\,j=2,4, \\[8pt]
|\A_{st}(\mathcal{Q}_{h}u^{m},v_{h}^{m})|& \leq C h^{k+1/2}\|u\|_{k+1}\vertiii{v_{h}^{m}}_{\mathcal{T}_{h}}.
\end{array}
\]
\end{lemma}
\noindent
{\bf Proof.}  By using $L^{2}$-projection error estimates given in \cite[Lemma 3.6]{lin2018weak}, we deduce that
\[
\begin{array}{ll}
|\mathcal{L}_{1}(u^{m},v_{h}^{m})|=\displaystyle \left| \sum_{T \in \mathcal{T}_{h}}\left(u^{m}-\mathcal{Q}_{h}u^{m},\s_{m}\cdot \nabla v_{h}^{m}\right)_{T}\right|
 \leq C h^{k+1}\|u^{m}\|_{k+1}\vertiii{v_{h}^{m}}_{\mathcal{T}_{h}}.
\end{array}
\]
Using Cauchy-Schwarz inequality and trace inequality (\ref{traceeq}), we have
\[
\begin{array}{ll}
|\mathcal{L}_{2}(u^{m},v_{h}^{m})|&=\displaystyle \left|\sum_{T \in \mathcal{T}_{h}}\left\langle u^{m}-\{\mathcal{Q}_{h}u^{m}\},\s_{m} \cdot \mathbf{n} (v_{h}^{m}-\{v_{h}^{m}\}) \right\rangle_{\partial T}\right|\\[16pt]
& \leq \displaystyle \sum_{T \in \mathcal{T}_{h}}\|u^{m}-\{\mathcal{Q}_{h}u^{m}\}\|_{\partial T}\||\s_{m} \cdot \mathbf{n}|^{1/2}(v_{h}^{m}-\{v_{h}^{m}\})\|_{\partial T} \\[16pt]
& \leq C h^{k+1/2}\|u^{m}\|_{k+1}\vertiii{v^{m}}_{\mathcal{T}_{h}}.
\end{array}
\]
In a similar way $\mathcal{L}_{3}(u^{m},v_{h}^{m})$ is bounded by
\[
\begin{array}{ll}
|\mathcal{L}_{3}(u^{m},v_{h}^{m})|&=\displaystyle  \sum_{T \in \mathcal{T}_{h}}\left| \left(\sigma_{t} (u^{m}-\mathcal{Q}_{h}u^{m})-\sigma_{s} \sum_{\ell =0}^{M}\w_{\ell}\Phi(\x,\s_{m}\cdot \s_{\ell})(u^{\ell}-\mathcal{Q}_{h}u^{\ell}) , v^{m}_{h}\right)_{T}\right|\\[16pt]
& \leq \displaystyle  C \sum_{T \in \mathcal{T}_{h}}\|u^{m}-\mathcal{Q}_{h}u^{m}\|_{ T}\|v^{m}_{h}\|_{ T}  \leq C h^{k+1}\|u^{m}\|_{k+1}\vertiii{v^{m}}_{\mathcal{T}_{h}},
\end{array}
\]
where  the assumption (\ref{quadasum}) is used. By the means of  trace inequality (\ref{traceeq}), 
\[
\begin{array}{ll}
|\mathcal{L}_{4}(u^{m},v_{h}^{m})| &\leq C \|u^{m}-\mathcal{Q}_{h}u^{m}\|_{\partial \Omega^{-}_{h}}\||\s_{m} \cdot \mathbf{n}|^{1/2}v_{h}^{m}\|_{\partial \Omega^{-}_{h}} \\ [8pt]
&\leq C h^{k+1/2}\|u^{m}\|_{k+1}\vertiii{v^{m}}_{\mathcal{T}_{h}}.
\end{array}
\]
At last, the stabilization term can be estimated by
\[
|\A_{st}(\mathcal{Q}_{h}u^{m},v^{m})|  \leq C h^{k+1/2}\|u^{m}\|_{k+1}\vertiii{v^{m}}_{\mathcal{T}_{h}}.
\]
This completes the proof.
\hfill \eop

\begin{theorem}\label{errest1}
Let $ u^{m} $ and 	$ u^{m}_{h} $ be the solutions to the problem (\ref{SRtde1}) and the discrete WG scheme (\ref{wg1}), respectively. Then, the error term
$\mathbf{e}_{\mathcal{Q}_{h}}\equiv \mathcal{Q}_{h}\bu-\bu_{h}$
satisfies:
\begin{equation}\label{numerrest}
\vertiii{\mathcal{Q}_{h}\bu-\bu_{h}} \leq C h^{k+1/2}\|\bu\|_{k+1}.
\end{equation}
\end{theorem}

\noindent
{\bf Proof.} Using the stability result from Lemma \ref{stab}, we have 
\begin{equation}\label{estpf1}
\alpha \vertiii{\mathcal{Q}_{h}\bu-\bu_{h}}^{2} \leq \A_{WG}(\mathcal{Q}_{h}\bu-\bu_{h},\mathcal{Q}_{h}\bu-\bu_{h}).
\end{equation}
By employing Lemma \ref{errest0}, one can obtain that
\begin{equation}\label{estpf2}
\begin{array}{ll}
\A_{WG}(\mathcal{Q}_{h}\bu-\bu_{h},\mathcal{Q}_{h}\bu-\bu_{h}) & = \mathbf{L}(\bu,\mathcal{Q}_{h}\bu-\bu_{h})\\[6pt]
& \leq C h^{k+1/2}\|\bu\|_{k+1}\vertiii{\mathcal{Q}_{h}\bu-\bu_{h}}\\[6pt]
& \leq C h^{2k+1}\|\bu\|^{2}_{k+1}+\dfrac{\alpha}{2}\vertiii{\mathcal{Q}_{h}\bu-\bu_{h}}^{2}.
\end{array}
\end{equation}
A straightforward computation provides the desired result (\ref{numerrest}).
\hfill \eop

\bigskip
The projection error is discussed here. 
By using the trace inequality (\ref{traceeq}) and result form  \cite[Lemma 3.6]{lin2018weak},  it can be easily verified that 
\begin{equation}\label{projerr}
\vertiii{u^{m}-\mathcal{Q}_{h}u^{m}}_{\mathcal{T}_{h}} \leq C h^{k+1/2}\|u^{m}\|_{k+1}.
\end{equation}

Next, we combine the above results to deduce the error estimate associated with weak Galerkin approximation in the following theorem.
\begin{theorem}\label{errest2}
Let  $\bu_{h}= \{u^{m}_{h}\} $ be the numerical approximation of $ \bu=\{u^{m}\}$. Then, the following estimates 
\begin{equation}\label{numerr}
\vertiii{\bu-\bu_{h}}\leq C h^{k+1/2}\|u^{m}\|_{k+1}.
\end{equation}
\end{theorem}

\noindent
{\bf Proof.} Combing the error estimates (\ref{numerrest}) and (\ref{projerr}), we obtain that
\begin{equation}\label{wdgerror}
\begin{array}{ll}
\vertiii{u^{m}-u^{m}_{h}}_{\mathcal{T}_{h}} &= \vertiii{u^{m}-\mathcal{Q}_{h}u^{m}}_{\mathcal{T}_{h}}+\vertiii{\mathcal{Q}_{h}u^{m}-u^{m}_{h}}_{\mathcal{T}_{h}} \\[8pt]
&\leq   C h^{k+1/2}\|u^{m}\|_{k+1}.
\end{array}
\end{equation}
Next, apply the norm definition from (\ref{normdef}) to derive (\ref{numerr}). This finalizes the proof.
\hfill \eop

\vskip 2mm Finally, we discuss the estimate for the fully discrete error $\mathbf{e}_{h}:=\{e_{h}^{m}\}_{m=0}^{M}$, where 
$e^{m}_{h} =  (u-u_{h}^{m}) $.
By adding  estimates from  Theorem \ref{domerrthm} and Theorem \ref{errest2}, the fully discrete   error term $\mathbf{e}_{h} \equiv \bu-\bu_{h}$  satisfies:

In the $2D$ model setting, 
\begin{equation}\label{finerror2d}
\vertiii{\mathbf{e}_{h}}=\left(\sum_{m=0}^{M}\w_{m}\vertiii{e^{m}_{h}}^{2}_{\mathcal{T}_{h}} \right)^{1/2}\leq C h^{k+1/2}\|u^{m}\|_{k+1}+ \mathcal{O}(h_{\theta}^2).
 \end{equation}
 In the $3D$ model setting, 
\begin{equation}\label{finerror3d}
\vertiii{\mathbf{e}_{h}}=\left(\sum_{m=0}^{M}\w_{m}\vertiii{e^{m}_{h}}^{2}_{\mathcal{T}_{h}}  \right)^{1/2}\leq C_{1} h^{k+1/2}\|u^{m}\|_{k+1}+ C_{2}\,n^{-r}\left(\int_{\Omega}\|u\|^{2}_{r,\Sc}\,\mbox{d}\x\right)^{1/2}.
\end{equation}

\section{Numerical discussion}\label{Sec4}

In this section, we consider a few test examples employing manufacturing methods to verify the theoretical estimations provided in the earlier sections.  For the theoretical validation, we focus the numerical experiments only on the two-dimensional models to lower the computational cost.   The numerical construction is followed from \cite{wang}, which is described briefly here. 

Firstly, we state that the spatial domain is  $\Omega = (0,1)^2$ and the mesh width in the angular direction (parameterized form) is  $h_{\theta} = \pi/10$ for all the numerical examples.
Let $ \mathcal{T}^{m}_{h}$ be the triangulation of $\Omega$ with respect to the angular direction $\s_{m}$. For example, the angular direction  $\s_{0}\equiv (\cos{0}, \sin{0})=(1,0)$ is associated with the triangulation $ \mathcal{T}^{0}_{h}$, which has the incoming boundary $\partial \Omega_{\s_0}^{-}=\{(0,y)|y=[0,1]\}$. In this way, we can fix the spatial triangulation  $ \mathcal{T}^{m}_{h}$ associated with the angular direction 
$\s_{m}$, then the following algorithm is used for the numerical simulation.

\begin{algorithm}[!ht]
\caption{Source iteration}\label{alg:si}
{\bf Set parameters} $tol, \delta$\\
{\bf Initialization} $\{u^{m,0}_{h}\} = \{0\},\,\,err = tol +1, \, i=0$\\
 \While{$err \geq tol$}{
  Set $\{u^{m,i+1}_{h}\}=\{u^{m,i}_{h}\}$\\
  Compute $\mathcal{K}_{d}u^{m,i}_{h} = \sum_{m=0}^{M}\w_{m}\Phi(\x,\s\cdot \s_{m})u^{m,i}_{h}$\\
  \For{$m=0$ to $M$}{
  Update $\{u^{m,i+1}_{h}\}$ on all elements $T \in \mathcal{T}^{m}_{h}$ by solving the discrete problem (\ref{wg1})
   }{
   Compute $err = \|\{u^{m,i+1}_{h}\}-\{u^{m,i}_{h}\}\|_{\Omega_{h}}$\\
   Set $i = i+1$
  }
 }
\end{algorithm}

We complete the aforementioned process in a single iteration step for all directions, terminate the iteration if a stopping condition (for below test problems, $tol = 10^{-3}$) is fulfilled, and take $u^{m,i}_{h}$ as $u^{m}_{h}$.   We have employed the software library deal.II \cite{wang2019deal} for the same computation.  All the numerical simulations were executed on an Intel Core i7-10510U CPU@1.80GHz machine with 8 GB of RAM.

\subsection{Numerical examples}
Let $\Omega = ( 0,1 )^{2}$. We consider the following  examples of the radiative transfer equation (\ref{SRte1}):
\begin{example}\label{numex1}
The H-G phase function (\ref{hgfun}) with $\eta = 0.5 $ and the exact solution is 
\[u(\x,\s)=\sin(\pi x_{1}) \sin(\pi x_{2})\]
\end{example}

\vskip 2mm

\begin{example}\label{numex2}
The  scattering phase function   $\Phi(t)=(2+t)/4 \pi$ and the exact solution is
\[u(\x,\s)=e^{-ax_{1}-bx_{2}}(1+c \cos(\theta)), \quad \theta \in [0, 2\pi],\] where $a,b=1/2$ and $c = 1/(1 + 6 \sigma_s)$. 

\vskip 2mm
For both test problems, the source functions $f$ satisfy the radiative
transfer equation (\ref{SRte1}). Furthermore,  the absorption and scattering coefficients value is $\sigma_{t}=2,\,\,\sigma_{s}=1/2$ for all the test problems.
\end{example}

For numerical validation, we discuss error $\mathbf{e}_{h}$ in the following norm 
\[
\vertiii{\mathbf{e}_{h}}_{1} = \left(\sum_{m =0}^{M}\w_{m}\sum_{T \in \mathcal{T}_{h}}\|e^{m}_{h}\|^{2}_{T} \right)^{1/2},
\]
and convergence order $eoc_{h}$ is given by
\[
eoc_{h} = \log_{2}\left(\dfrac{\vertiii{\mathbf{e}_{h}}_{1}}{\vertiii{\mathbf{e}_{h/2}}_{1}}\right).
\]

In spatial discretization, the coarse mesh (at first level)  consists of four squares. Then, each square is equally decomposed into $4$ squares, to get the next fine mesh. For the numerical experiments, $Q_1$, the bilinear ﬁnite elements on quadrilaterals, and $Q_2$, the biquadratic ﬁnite elements on quadrilaterals are chosen and error estimates with convergence orders are derived.

In table \ref{tab1}, the error estimate, for example, \ref{numex1} is presented whereas the error estimate for Example \ref{numex2} in table \ref{tab2}. Also, the order of convergence ($eoc_h$) is calculated for the error estimates. From tables \ref{tab1}, \ref{tab2}, one can verify  the theoretical error estimate of  $\mathcal{O}(h^{k+1/2})$ with $Q_1$ and $Q_2$-finite elements. In particular, $eoc_h$ is close to $2.5$ rather than $3$ for  $Q_2$-finite elements. This reduction happens due to the composite trapezoidal formula which is of second order.

\begin{table}[!ht]
	\renewcommand{\arraystretch}{1}
	\caption{\label{tab1}  \it{Discretization error $e_h$ and order of convergence $eoc_h$  for  Example \ref{numex1}. }}
	{\centering
		\begin{tabular}{c|| c c ||c c }	
		\hline
			$1/h$  & $Q_1$-element    &  & $Q_2$-element   &   \\
			\hspace{2cm} &   $\vertiii{\mathbf{e}_{h}}_{1}$  & $eoc_h$ &  $\vertiii{\mathbf{e}_{h}}_{1}$  & $eoc_h$\\
			\hline &&&&\\[-12pt]
			8  & 8.5643e-02 & -- & 6.4527e-02& -- \\[6pt]
			16	 &2.7626e-02  & 1.63  & 1.2936e-02& 2.32 \\[6pt]
			32 &  8.4227e-03 &  1.71  &2.4672e-03& 2.39 \\[6pt]
			64	&  2.4413e-03  & 1.78  &4.4057e-04&2.48\\[6pt] 
			128	&  6.7628e-04  & 1.85  &7.4519e-05&2.56\\  
			\hline
		\end{tabular}
		\par}
\end{table}

\begin{table}[!ht]
	\renewcommand{\arraystretch}{1}
	\caption{\label{tab2}  \it{Discretization error $e_h$ and order of convergence $eoc_h$   for  Example \ref{numex2}. }}
	{\centering
		\begin{tabular}{c|| c c ||c c }	
		\hline
			$1/h$  & $Q_1$-element    &  & $Q_2$-element   &   \\
			\hspace{2cm} &   $\vertiii{\mathbf{e}_{h}}_{1}$  & $eoc_h$ &  $\vertiii{\mathbf{e}_{h}}_{1}$  & $eoc_h$\\
			\hline &&&&\\[-12pt]
			8  & 7.5954e-02 & -- & 5.7826e-02& -- \\[6pt]
			16 &2.5150e-02  & 1.59  & 1.1892e-02& 2.28 \\[6pt]
			32	&  8.0352e-03 &  1.64  &2.3037-03& 2.36 \\[6pt]
			64 &  2.4497e-03  & 1.72  &4.2676e-04& 2.43\\[6pt]
			128	&  6.9794e-03  & 1.81  &7.4895e-05& 2.51\\ 
			\hline
		\end{tabular}
		\par}
\end{table}

\subsection{Comparison with DG methods}
In this part of our article, we discuss a detailed comparison of our proposed DOWG method with the existing DG-type methods such as the  discrete-ordinate discontinuous-Galerkin (DODG) method,  the  discrete-ordinate discontinuous-streamline
diffusion (DODSD method) presented in \cite{han,wang}. A comparative analysis of DG and WG methods for the elliptic boundary value problem is discussed in \cite{lin2015comparative}.

\subsubsection{DODG method for RTE model (\ref{SRte1})}
DODG method for the model problem (\ref{SRte1}) is given by
\begin{equation}\label{dg1}
\left\{
\begin{array}{ll}
\mbox{Find}\,\, \bu_{h} \in \mathbf{V}_{h} \,\,\mbox{such that}\\[6pt]
\A_{DG}(\bu_{h},\bv_{h})=\F(\bv_{h}), \,\,\forall \, \bv_{h} \in \mathbf{V}_{h},
\end{array}
\right.
\end{equation}
where 
\[
\A_{DG}(\bu_{h},\bv_{h})=\sum_{m=0}^{M}\w_{m} \big(\A_{cdg}(u_{h}^{m},v^{m}_{h})+\A_{stdg}(u_{h}^{m},v_{h}^{m})\big),\quad \F(\bv_{h})=\sum_{m=0}^{M}\w_{m} \F^{\s}(v_{h}^{m}),
\]
with $\F^{\s}(v_{h}^{m})$ is given as in (\ref{SRtfe1}) and 
\[
\displaystyle \A_{cdg}(u_{h}^{m},v^{m}_{h}) =  \bigg(\s_{m} \cdot \nabla_{w} u_{h}^{m}+\sigma_{t} u_{h}^{m}-\sigma_{s} \sum_{\ell=0}^{M}\w_{\ell}\Phi(\x,\s_{m}\cdot \s_{\ell})u_{h}^{\ell},v_{h}^{m}\bigg)-\big\langle \s_{m} \cdot \mathbf{n} \,\hat{u}_{h}^{m},v_{h}^{m}\big\rangle_{\partial \Omega_{h}},
\]
\[
\A_{stdg}(u_{h}^{m},v^{m}_{h}) = c_{p} \sum_{e \in \mathcal{E}^{0}_{h}}\big\langle[u_{h}^{m}], [v_{h}^{m}]\big\rangle_{e}.
\]
Here, $c_p > 0$ is  a penalty parameter and  $\hat{u}_{h}^{m}$ is the numerical trace, which is  given by
\[
\hat{u}_{h}^{m} = 
\left\{
\begin{array}{ll}
     u_{in}^{m} &  \mbox{in} \,\,\Gamma^{-}_{m}  \\
     \lim_{\epsilon \rightarrow 0+}u_{h}^{m}(\x-\epsilon \s_{m})& \mbox{otherwise}.
\end{array}
\right.
\]

\subsubsection{DODSD method for RTE model (\ref{SRte1})}
DODSD method for the model problem (\ref{SRte1}) is given by
\begin{equation}\label{sd1}
\left\{
\begin{array}{ll}
\mbox{Find}\,\, \bu_{h} \in \mathbf{V}_{h} \,\,\mbox{such that}\\[6pt]
\A_{DSD}(\bu_{h},\bv_{h})=\F(\bv_{h}), \,\,\forall \, \bv_{h} \in \mathbf{V}_{h},
\end{array}
\right.
\end{equation}
where 
\[
\A_{DSD}(\bu_{h},\bv_{h})=\sum_{m=0}^{M}\w_{m} \big(\A_{cdsd}(u_{h}^{m},v^{m}_{h})+\A_{stdsd}(u_{h}^{m},v_{h}^{m})\big),\quad \F(\bv_{h})=\sum_{m=0}^{M}\w_{m} \F^{\s}(v_{h}^{m}+\delta \s_{m} \cdot \nabla v_{h}^{m}),
\]
with
\[
\displaystyle \A_{cdsd}(u_{h}^{m},v^{m}_{h}) =  \bigg(\s_{m} \cdot \nabla_{w} u_{h}^{m}+\sigma_{t} u_{h}^{m}-\sigma_{s} \sum_{\ell=0}^{M}\w_{\ell}\Phi(\x,\s_{m}\cdot \s_{\ell})u_{h}^{\ell},v_{h}^{m}+\delta \s_{m} \cdot \nabla_{w} v_{h}^{m}\bigg)+\sum_{T \in \mathcal{T}_{h}}\big\langle[u^{m}_{h}], v^{m,+}_{h}|\s_{m}\cdot \mathbf{n}|\big\rangle_{\partial^{-}T},
\]
and $u^{m}_{h} = u^{m}_{in}$ on $\Gamma^{-}_{m}$. The stabilization parameter $\delta$ satisfies  $\delta= c h$ with some $c > 0$.

Note that, the numerical solutions $u_{h}$ of both DG-type methods satisfy the error estimates (\ref{finerror2d}) and (\ref{finerror3d}).

Now, we discuss the comparison of these methods with the proposed method (\ref{wg1})
As, we can see from both numerical schemes (\ref{dg1}) and (\ref{sd1}), these methods are parameter-dependent. In contrast, the proposed method (\ref{wg1}) is parameter-free, one of the major advantages of our numerical method. In the numerical section of \cite{han}, the authors discussed the impact of parameter $c_p$ on the accuracy of the numerical method.

To compare the accuracy of numerical solution with the exact solution of DG methods versus WG methods, we discuss the error estimate   $\vertiii{\mathbf{e}_{h}}_{1}$  for the $Q_1$ finite elements. The penalty parameter $c_p = 0.1$ and stabilization parameter $\delta = h$ is taken for numerical methods (\ref{dg1}) and (\ref{sd1}), respectively.


\begin{table}[!ht]
\caption{\label{tab3}  \it{Comparison of error $\vertiii{\mathbf{e}_{h}}_{1}$  for  Example \ref{numex1}.}}
    \centering
    \begin{tabular}{c||cc||cc||cc}
    \hline
        $1/h$ & WDG method && DG method  && DSD method &\\
        \hline &&&&&&\\[-12pt]
        8  & 8.5643e-02 & --  & 8.9306e-02 & -- & 8.6366e-02&-- \\
        16  & 2.7626e-02  & 1.63    &3.1009e-02 &1.52 &2.8981e-02& 1.57 \\
        32  & 8.4227e-03 &  1.71   & 9.9387e-03 & 1.64&9.1714e-03&1.66\\
        64  & 2.4413e-03  & 1.78   & 3.0301e-03 & 1.71 &2.7459e-03& 1.74\\
        128  & 6.7628e-04  & 1.85  & 8.7072e-04 & 1.80&7.8232e-04&1.81\\
        \hline
    \end{tabular}
    \label{compare_err_ex1}
\end{table}


\begin{table}[!ht]
\caption{\label{tab4}  \it{Comparison of error $\vertiii{\mathbf{e}_{h}}_{1}$  for  Example \ref{numex2}.}}
    \centering
    \begin{tabular}{c||cc||cc||cc}
    \hline
        $1/h$ & WDG method && DG method  && DSD method &\\
        \hline &&&&&&\\[-12pt]
        8  & 7.5954e-02 & --  & 7.7445e-02 & & 7.6240e-02& --\\
        16  & 2.5150e-02  & 1.59     &2.6265e-02 & 1.56&2.5583e-02& 1.57\\
        32  &  8.0352e-03 &  1.64  & 8.4862e-03 & 1.63&8.2528e-03&1.63\\
        64  & 2.4497e-03  & 1.72   & 2.6119e-03 & 1.70&2.5199e-03&1.71 \\
        128  & 6.9794e-04  & 1.81  & 7.6149e-04 & 1.78&7.3042e-04&1.78\\
        \hline
    \end{tabular}
    \label{compare_err_ex2}
\end{table}

We have presented error estimates in tables \ref{tab3} and \ref{tab4} for all three numerical methods, one can see that our proposed  WG method performs better than other DG methods. The convergence rate of our method is also slightly better than both DG-type methods. 
We conclude this subsection with a few remarks that the DOWG method is a viable alternative to DODG/DODSD methods. Compared to DG-type methods, it is easier to use and performs better and no need for penalty/stabilization parameters. However, a complete discussion of this topic is beyond the scope of this paper.

\section{Concluding remarks}\label{Sec5}
In this article, we have proposed a weak Galerkin finite element method for the numerical solution of the radiative transfer equation. The discrete-ordinate weak Galerkin method has been used to approximate the angular and spatial variables. For the resulting numerical method, the stability result is established. In the specific norm,   the error estimate is devised for the numerical solutions under the suitable regularity assumptions. In the end, the numerical experiments are presented for the justification of theoretical findings.

\section*{Acknowledgments}
The author is grateful to the National Board of Higher Mathematics, Department of Atomic Energy (DAE), India for granting postdoctoral fellowships during his tenure at the Indian Institute of Science, Bangalore. 

\bibliography{MKS_WG_RTE}







\end{document}